\documentclass[11pt]{amsart}
\usepackage{geometry}                % See geometry.pdf to learn the layout options. There are lots.
\usepackage{graphicx}
\usepackage{amssymb}
\usepackage{amsthm}
\usepackage{mathrsfs}
\usepackage{epstopdf}
\usepackage{MnSymbol}
\usepackage{stmaryrd}
\usepackage{verbatim}
\usepackage[all]{xy}
\title[A stable but fibrewise nonstable bundle on the twistor space]{An example of a stable but fibrewise nonstable bundle on the twistor space of \\ a hyperk\"ahler manifold}
\author{Artour Tomberg}
\thanks{The study has been funded by the Russian Academic Excellence Project '5-100'}
\address{Faculty of Mathematics, National Research University Higher School of Economics, 6 Usacheva St., Moscow, Russia, 119048}
\address{Department of Mathematics, Western University, Middlesex College, London, Ontario, Canada, N6A 5B7}
\email{artour@tomberg.com}

\begin{document}

\begin{abstract}
We construct an explicit example of a stable bundle on the twistor space $\mathrm{Tw}(M)$ of a hyperk\"ahler manifold $M$ whose restrictions to all the fibres of the natural twistor projection $\pi : \mathrm{Tw}(M) \to \mathbb{CP}^1$ are nonstable. We also describe the relationship between bundles on $\mathrm{Tw}(M)$ that do not have subsheaves of strictly lower rank and bundles that stably restrict to the fibres of $\pi$, and announce a result whose proof will appear in a forthcoming paper.
\end{abstract}

\maketitle

\tableofcontents

\section{Introduction} A \emph{hyperk\"ahler manifold} is a smooth manifold $M$ together with a triple of integrable almost complex structures $I, J, K : TM \to TM$ satisfying the quaternionic relations $I^2 = J^2 = K^2 = -1$, $IJ = -JI = K$, and a Riemannian metric $g$ which is K\"ahler with respect to the structures $I, J, K$. An example of a hyperk\"ahler manifold is a \emph{K3 surface}, that is, a compact simply connected complex surface with trivial canonical bundle. It admits a hyperk\"ahler metric as a consequence of the Calabi-Yau theorem \cite{yau}. 

It's not hard to see that a hyperk\"ahler manifold $M$ admits a whole family of \emph{induced complex structures}, which topologically looks like a 2-sphere:
\[
S^2  = \left\{aI + bJ + cK : a^2 + b^2 + c^2 = 1 \right\}.
\]
Identifying $S^2$ with $\mathbb{CP}^1$, we define the \emph{twistor space} of $M$ as the topological Cartesian product $\mathrm{Tw}(M) := M \times S^2 = M \times \mathbb{CP}^1$. We think of $\mathrm{Tw}(M)$ as parametrizing the induced complex structures at points of $M$. In the context of the twistor space, the initial structures $I, J, K$ don't play any vital role, and henceforth we will denote by $I \in \mathbb{CP}^1$ an arbitrary induced complex structure, while $M_I$ will denote the corresponding complex manifold. Note that for any $I \in \mathbb{CP}^1$, $g$ is a K\"ahler metric on $M_I$ .

The twistor space $\mathrm{Tw}(M)$ admits a natural integrable almost complex structure \cite{salamon}.
With respect to this structure, the projection onto the second coordinate $\pi : \mathrm{Tw}(M) \to \mathbb{CP}^1$ is holomorphic, and the fibres of $\pi$ over points $I \in \mathbb{CP}^1$ correspond to the complex manifolds $M_I$. We observe that the projection onto the first coordinate of the twistor space $\sigma : \mathrm{Tw}(M) \to M$ is not holomorphic with respect to any of the induced complex structures on $M$, since $\mathrm{Tw}(M)$ is not a product of $M$ and $\mathbb{CP}^1$ as complex manifolds, but only as topological manifolds. $\mathrm{Tw}(M)$ also admits a natural Hermitian metric satisfying the \emph{balancedness} condition $d\left(\omega^{n-1}\right) = 0$, where $\omega$ is the Hermitian form of this metric and $n$ is the complex dimension of $\mathrm{Tw}(M)$ \cite{kaled-verbit}.

Let $E$ be a holomorphic vector bundle on the twistor space $\mathrm{Tw}(M)$ of a compact hyperk\"ahler manifold $M$. We define the \emph{degree} of $E$ by 
\begin{equation} \label{degree}
\deg (E) =\int_{\mathrm{Tw}(M)} c_1(E) \wedge \omega^{n-1},
\end{equation}
where by $c_1(E)$ we denote any representative of the first Chern class of $E$ in $H^2(M, \mathbb{R})$. We say that $E$ is \emph{stable} if for any subsheaf $\mathcal{F} \subset E$ satisfying $0 < \mathrm{rk}(\mathcal{F}) < \mathrm{rk}(E)$, we have strict inequality
\[
\frac{\deg(\mathcal{F})}{\mathrm{rk}(\mathcal{F})} < \frac{\deg(E)}{\mathrm{rk}(E)},
\]
where the degree of $\mathcal{F}$ is defined as the degree of its determinant line bundle: $\deg(\mathcal{F}) := \deg(\det \mathcal{F})$. The bundle $E$ is called \emph{irreducible} if it does not have any nonzero subsheaves of lower rank.

Observe that the value of the integral \eqref{degree} in the definition of degree does not depend on the choice of the representative of the Chern class $c_1(E)$, since the Hermitian form $\omega$ on $\mathrm{Tw}(M)$ satisfies the balancedness condition $d\left(\omega^{n-1}\right) = 0$, as we noted earlier. For every $I \in \mathbb{CP}^1$, we can similarly define the degree of bundles on the K\"ahler manifold $M_I$, since any K\"ahler metric is a priori balanced. In this way, the notion of degree makes sense both for bundles $E$ on the twistor space $\mathrm{Tw}(M)$ and their restrictions $E_I$ to the fibres $\pi^{-1}(I) = M_I$ of the twistor projection $\pi: \mathrm{Tw}(M) \to \mathbb{CP}^1$.

In the paper \cite{kaled-verbit}, Kaledin and Verbitsky prove the following result, among other things.

\theoremstyle{plain}
\newtheorem{forward_weak}{Proposition}
\begin{forward_weak} \label{init-version}
Let $M$ be a compact hyperk\"ahler manifold and $E$ a holomorphic vector bundle on the twistor space $\mathrm{Tw}(M)$. If $E$ stably restricts to the generic fibre of the holomorphic twistor projection $\pi: \mathrm{Tw}(M) \to \mathbb{CP}^1$ (in the sense of the Zariski topology on $\mathbb{CP}^1$), then it is stable as a bundle on $\mathrm{Tw}(M)$ as well.
\end{forward_weak}

In the present short note, we will show that the converse statement does not not hold in general. More precisely, we will construct an explicit example of a stable holomorphic bundle $E$ of rank 2 on the twistor space $\mathrm{Tw}(M)$ of a K3 surface $M$, all of whose restrictions to the fibres of the projection $\pi: \mathrm{Tw}(M) \to \mathbb{CP}^1$ are nonstable. We will also formulate a stronger version of the above result (Theorem \ref{forward}) and briefly discuss the converse to this stronger statement.

\section{An example of a stable but fibrewise nonstable bundle on $\mathrm{Tw}(M)$}

Let $M$ be an algebraic K3 surface with Picard number $\rho(M) \ge 2$ (for basic properties of K3 surfaces, as well as terminology and important theorems in complex geometry that we use below, see e.g. \cite{griffiths-harris}). The degree of any bundle on $M$ is an integer, hence for line bundles we have a homomorphism of groups $\deg : \mathrm{Pic}(M) \to \mathbb{Z}$ with a nonzero kernel. We can choose an element $L$ of this kernel in such a way that the inequality $h^1(M, L^*) \ne 0$ holds. Indeed, by the Riemann-Roch formula for K3 surfaces,
\[
h^0(M, L) - h^1(M, L) + h^2(M, L) = \frac{c_1(L)^2}{2} + 2.
\]
Let $L$ be a nontrivial holomorphic line bundle of degree zero. $L$ does not have any nonzero sections, since such a section would give an effective divisor of a strictly positive degree, by the Poincar\'e-Lelong formula. Hence $h^0(M, L) = 0$, and similarly $h^2(M, L) = h^0(M, L^*) = 0$, where we use Serre duality. Thus,
\[
h^1(M, L) = - \frac{c_1(L)^2}{2} - 2,
\]
and by the Hodge index theorem, $c_1(L)^2 < 0$. Replacing $L$ by its multiple, if necessary, we have $h^1(M, L) \ne 0$, and thus $h^1(M, L^*) \ne 0$.

Since $\deg L = 0$, one can show (see \cite{verbit1}, Theorem 2.4) that $L$ is \emph{hyperholomorphic}, that is, admits a Hermitian metric with Chern connection $\nabla$ whose curvature is a $(1,1)$-form with respect to every induced complex structure on $M$. Clearly, this means that for every $I \in \mathbb{CP}^1$, $\nabla$ endows $L$ with the structure of a holomorphic line bundle over $M_I$, which we will denote by $L_I$. Moreover, taking the pullback of $(L, \nabla)$ along the projection onto the first coordinate of the twistor space $\sigma : \mathrm{Tw}(M) \to M$, we get a line bundle $\sigma^*L$ on $\mathrm{Tw}(M)$ with holomorphic structure $(\sigma^*\nabla)^{0,1}$. The restriction of the holomorphic line bundle $\sigma^* L$ to the fibre $\pi^{-1}(I) = M_I$ of the twistor projection $\pi: \mathrm{Tw}(M) \to \mathbb{CP}^1$ is precisely $L_I$. We will denote the initial complex manifold structure on our K3 surface (which corresponds to one of the $I \in \mathbb{CP}^1$) simply by $M$, while the initial holomorphic structure on our line bundle will be denoted by $L$, rather than $L_I$; this should not cause any confusion.

The higher direct images of $\sigma^* L^*$ with respect to the projection $\pi : \mathrm{Tw}(M) \to \mathbb{CP}^1$ are as follows (see \cite{verbit3}, Proposition 6.3):
\[
R^i\pi_* \left( \sigma^* L^* \right) \cong \mathcal{O}_{\mathbb{CP}^1}(i) \otimes_{\mathbb{C}} H^i(M, L^*).
\]
Let us denote $\mathcal{O}_{\mathrm{Tw}(M)}(-1) := \pi^*\left( \mathcal{O}_{\mathbb{CP}^1}(-1) \right)$. Applying the projection formula and the above,
\[
R^1\pi_*\left[\sigma^*L^* \otimes \mathcal{O}_{\mathrm{Tw}(M)}(-1)\right] \cong \left[R^1\pi_*(\sigma^* L^*)\right] \otimes \mathcal{O}_{\mathbb{CP}^1}(-1) \cong \mathcal{O}_{\mathbb{CP}^1} \otimes_{\mathbb{C}} H^1(M, L^*).
\]
Thus,
\[
\mathrm{Ext}^1(\sigma^* L, \mathcal{O}_{\mathrm{Tw}(M)}(-1)) \cong H^1(\sigma^* L^* \otimes \mathcal{O}_{\mathrm{Tw}(M)}(-1)) =
\]
\[
= \Gamma\left(\mathbb{CP}^1, R^1\pi_*\left[\sigma^*L^* \otimes \mathcal{O}_{\mathrm{Tw}(M)}(-1)\right]\right) \cong H^1(M, L^*).
\]
This is nonzero by construction, so we can choose a nonzero element in $H^1(M, L^*) = \mathrm{Ext}^1(\sigma^* L, \mathcal{O}_{\mathrm{Tw}(M)}(-1))$ which corresponds to some extension
\begin{equation} \label{extension}
0 \longrightarrow \mathcal{O}_{\mathrm{Tw}(M)}(-1) \longrightarrow E \longrightarrow \sigma^* L \longrightarrow 0.
\end{equation}
Observe that the restriction of this short exact sequence to any fibre $\pi^{-1}(I) = M_I$ of the holomorphic twistor projection $\pi: \mathrm{Tw}(M) \to \mathbb{CP}^1$ is $0 \to \mathcal{O}_{M_I} \to E_I \to L_I \to 0$. Since $\mathcal{O}_{M_I}$ and $L_I$ both have degree zero, $\deg E_I = \deg \mathcal{O}_{M_I} + \deg L_I$ is zero as well. Therefore, the morphism $\mathcal{O}_{M_I} \to E_I$ gives a destabilizing subsheaf of $E_I$, proving that $E_I$ is nonstable as a bundle on $M_I$.

We will now show that $E$ is stable as a bundle on $\mathrm{Tw}(M)$. One can show (see \cite{kaled-verbit}, Lemma 6.2) that the degree of any bundle on $\mathrm{Tw}(M)$ is equal to the degree of its restriction to any \emph{horizontal twistor line} $\left\{ m \right\} \times \mathbb{CP}^1 \subset \mathrm{Tw}(M)$, where $m \in M$. Clearly, the restriction of the exact sequence \eqref{extension} to any such line has the form
\[ \label{restriction}
0 \longrightarrow \mathcal{O}_{\mathbb{CP}^1}(-1) \longrightarrow \left.E\right|_{\left\{ m \right\} \times \mathbb{CP}^1} \longrightarrow \mathcal{O}_{\mathbb{CP}^1} \longrightarrow 0.
\]
This implies that $\deg E = \deg \mathcal{O}_{\mathbb{CP}^1}(-1) + \deg \mathcal{O}_{\mathbb{CP}^1} = -1 + 0 = -1.$ Moreover, since $\mathrm{Ext}^1(\mathcal{O}_{\mathbb{CP}^1}, \mathcal{O}_{\mathbb{CP}^1}(-1)) = H^1(\mathcal{O}_{\mathbb{CP}^1}(-1)) = 0$, we have $\left.E\right|_{\left\{m\right\} \times \mathbb{CP}^1} \cong \mathcal{O}_{\mathbb{CP}^1} \oplus \mathcal{O}_{\mathbb{CP}^1}(-1)$. This means that any potential destabilizing line subsheaf $\tilde{L} \hookrightarrow E$, that is, one which satisfies the inequality
\[
\deg(\tilde{L}) = \frac{\deg(\tilde{L})}{\mathrm{rk}(\tilde{L})} \ge \frac{\deg(E)}{\mathrm{rk}(E)} = -\frac{1}{2},
\]
should have degree 0. Let $\tilde{L} \hookrightarrow E$ be such a subsheaf. In the diagram
\begin{equation} \label{diagramma}
\xymatrix{&& \tilde{L} \ar[d] \ar[dr]^\theta && \\ 0 \ar[r] & \mathcal{O}_{\mathrm{Tw}(M)}(-1) \ar[r] & E \ar[r] & \sigma^* L \ar[r] & 0,}
\end{equation}
if the morphism $\theta : \tilde{L} \to \sigma^* L$ is zero, then by the exactness of the bottom row, there exists a lifting of the sheaf monomorhphism $\tilde{L} \hookrightarrow E$ to a monomorphism $\tilde{L} \hookrightarrow \mathcal{O}_{\mathrm{Tw}(M)}(-1)$. However, such a monomorphism cannot exist, since restricting any morphism $\tilde{L} \to \mathcal{O}_{\mathrm{Tw}(M)}(-1)$ to any horizontal twistor line $\left\{ m \right\} \times \mathbb{CP}^1$, we get $\mathcal{O}_{\mathbb{CP}^1} \to \mathcal{O}_{\mathbb{CP}^1}(-1)$, which is zero. On the other hand, if the morphism $\theta : \tilde{L} \to \sigma^* L$ is nonzero, it must be an isomorphism, since its restriction to any horizontal twistor line $\left\{ m \right\} \times \mathbb{CP}^1$ has the form $\mathcal{O}_{\mathbb{CP}^1} \to \mathcal{O}_{\mathbb{CP}^1}$, and any such nonzero morphism is an isomorphism. But if $\theta$ is an isomorphism, the diagram \eqref{diagramma} gives a splitting of the short exact sequence \eqref{extension}, which contradicts our choice of $E$. We have proved that such a destabilizing subsheaf $\tilde{L} \hookrightarrow E$ cannot exist, hence $E$ is stable.

\section{Irreducible bundles and fibrewise stability} The bundle $E$ on $\mathrm{Tw}(M)$ that we constructed in the previous section gives a counterexample to the converse of Proposition \ref{init-version}. However, looking at the proof of Lemma 7.3 in \cite{kaled-verbit}, it's not hard to see that Proposition \ref{init-version} can be made stronger in the following way.

\theoremstyle{plain}
\newtheorem{forward_strong}{Theorem}
\begin{forward_strong} \label{forward}
Let $M$ be a compact hyperk\"ahler manifold and $E$ a holomorphic vector bundle on the twistor space $\mathrm{Tw}(M)$. If $E$ stably restricts to the generic fibre of the holomorphic twistor projection $\pi: \mathrm{Tw}(M) \to \mathbb{CP}^1$, then it is irreducible as a bundle on $\mathrm{Tw}(M)$.
\end{forward_strong}

We believe that the converse to this stronger version of the statement is in fact true. At the present time, the following partial result is known.

\theoremstyle{plain}
\newtheorem{converse}[forward_strong]{Theorem}
\begin{converse} \label{converse}
Let $M$ be a compact simply connected hyperk\"ahler manifold and $E$ a holomorphic vector bundle on the twistor space $\mathrm{Tw}(M)$. If $E$ is irreducible, then it stably restricts to the generic fibre of the holomorphic twistor projection $\pi: \mathrm{Tw}(M) \to \mathbb{CP}^1$, provided that the rank of $E$ is equal to 2 or 3, or if its restriction $E_I$ to the generic fibre $\pi^{-1}(I) = M_I$ of $\pi$ is a simple bundle, in the sense that $\mathrm{Hom}(E_I, E_I) = \mathbb{C}$.
\end{converse}

The proof of this result will be given in a forthcoming paper. In the present short note, we will content ourselves with only a brief survey of the proof. We will make use of the following construction. For an arbitrary holomorphic vector bundle $E$ on $\mathrm{Tw}(M)$ and any $0 < s < \mathrm{rk}(E)$, we define the \emph{cone of exterior monomials} $C^s(E) \subseteq \Lambda^s E$ in the following way: at a point $x \in \mathrm{Tw}(M)$, $C^s(E)_x$ consists of the elements of $\Lambda^s E_x$ of the form $v_1 \wedge \ldots \wedge v_s$, where $v_1, \ldots, v_s \in E_x$. If $\mathcal{F} \hookrightarrow E$ is a subsheaf of rank $s$, it's not hard to verify that the image of $L = \det(\mathcal{F}) = \left(\Lambda^s \mathcal{F} \right)^{**} \hookrightarrow \left(\Lambda^s E\right)^{**} = \Lambda^s E$ lies in $C^s(E)$. At points where $\mathcal{F}$ is a subbundle of $E$, the line $L_x \subseteq \Lambda^s E_x$ is obtained from $\mathcal{F}_x \subseteq E_x$ by virtue of the Pl\"ucker embedding. On the other hand, one can show (see \cite{teleman}, subsection 2.2) that starting from a line subsheaf $L \subseteq \Lambda^s E$ with image in $C^s(E)$, one can recover a subsheaf $\mathcal{F} \subseteq E$ of rank $s$. Obviously, the above also holds for bundles on the fibres $\pi^{-1}(I) = M_I$ of the projection $\pi: \mathrm{Tw}(M) \to \mathbb{CP}^1$, and more generally on any complex manifold.

Let $M$ be a compact simply connected hyperk\"ahler manifold. It can be shown that for a bundle $E$ on $\mathrm{Tw}(M)$, viewed as a family of bundles on the fibres of the projection $\pi: \mathrm{Tw}(M) \to \mathbb{CP}^1$, stability is a Zariski open condition on $\mathbb{CP}^1$. In other words, the set of $I \in \mathbb{CP}^1$ for which the restriction $E_I$ is stable is Zariski open in $\mathbb{CP}^1$. The proof of this statement is essentially an adaptation of the argument from the proof of Theorem 1.3 in \cite{teleman}, where an analogous statement is shown for a projection $X \times Y \to X$, where $X \times Y$ is the product of complex manifolds $X$ and $Y$ satisfying certain properties.

Let $E$ be an irreducible bundle of rank $r$ on $\mathrm{Tw}(M)$. Arguing by contradiction, we assume that $E_I$ is nonstable as a bundle on $\pi^{-1}(I) = M_I$ for infinitely many $I \in \mathbb{CP}^1$. By the previous paragraph, it follows from this that $E_I$ is nonstable for all $I$. There exists $0 < s < r$ such that for every $I \in \mathbb{CP}^1$ there are destabilizing subsheaves $\mathcal{F}_I \hookrightarrow E_I$ of rank $s$, which correspond to line subsheaves $L_I \hookrightarrow \Lambda^s E_I$ with image in $C^s(E_I)$. Moreover, one can show that for some choice of $\mathcal{F}_I$, all the line bundles $L_I$ on $\pi^{-1}(I) = M_I$ are restrictions of a single line bundle $L$ on $\mathrm{Tw}(M)$.

Our goal consists in ``gluing'' these subsheaves $L_I \hookrightarrow \Lambda^s E_I$ over $M_I$ into a global subsheaf of $\Lambda^s E$ over $\mathrm{Tw}(M)$ with image in $C^s(E)$, from which we can recover a subsheaf of $E$ of rank $s$ over $\mathrm{Tw}(M)$ and get a contradiction to the irreducibility of $E$. Consider the vector bundle $\mathrm{Hom}(L, \Lambda^s E)$ on $\mathrm{Tw}(M)$. Its direct image $\pi_*\left(\mathrm{Hom}(L, \Lambda^s E)\right)$ along the twistor projection $\pi: \mathrm{Tw}(M) \to \mathbb{CP}^1$ is a vector bundle on $\mathbb{CP}^1$, and by Grauert's theorem (see \cite{grauert-remmert}, Theorem 10.5.5), $\pi_*\left(\mathrm{Hom}(\sigma^* L, \Lambda^s E)\right)_I = \mathrm{Hom}(L_I, \Lambda^s E_I)$ at all points $I \in \mathbb{CP}^1$, except perhaps finitely many. The existence of a subsheaf of $\Lambda^s E$ over $\mathrm{Tw}(M)$ with image in $C^s(E)$ will follow from the existence of the following algebraic morphism over $\mathbb{CP}^1$:
\[
\xymatrix{Y = \left\{\left[L_I \hookrightarrow C_s(E_I) \subseteq {{}\Lambda}^s E_I\right]\right\}_ {I \in \mathbb{CP}^1} \ar@{^{(}->}[r] \ar[d] & \mathbb{P}\left(\pi_*\left(\mathrm{Hom}(L, \Lambda^s E)\right)\right) \ar[dl] \\ \mathbb{CP}^1 &}
\]
If $r = $ 2 or 3, such a section always exists, since in this case $C^s(E) = \Lambda^s E$ for any $s$. If $r > 3$, the existence of a section of the morphism $Y \to \mathbb{CP}^1$ is not guaranteed, but such a section always exists over some ramified cover $f : X \to \mathbb{CP}^1$. Taking the fibred product
\[
\xymatrix{Z \ar[r]^-{\varphi} \ar[d]_\rho & \mathrm{Tw}(M) \ar[d]^\pi \\ X \ar[r]_f & \mathbb{CP}^1,}
\]
one can then proceed to construct a subsheaf $\mathcal{F} \subseteq \varphi^* E$ of rank $s$ over $Z$. If we assume that the restriction of $E$ to the generic fibre of $\pi : \mathrm{Tw}(M) \to \mathbb{CP}^1$ is a simple bundle, then after some work it can be shown that the irreducibility of $E$ on $\mathrm{Tw}(M)$ implies that $\varphi^* E$ is irreducible on $Z$, which leads to a contradiction.

The author would like to thank Misha Verbitsky and Dmitry Kaledin for their help in the preparation of the present note.

\end{document}